# *M|G|∞* Queue Parameters Values Approximation Through a Markov Renewal Process

**Manuel Alberto M. Ferreira**

Instituto Universitário de Lisboa (ISCTE – IUL), ISTAR- IUL, Lisboa, Portugal

**Abstract:** Some *M/G/∞* queue systems parameters values approximations, obtained through the consideration of an adequate Markov renewal process, are presented, and studied.

**Keywords:** *M/G/∞;* Markov renewal process; queue.

## 1. Introduction

In the *M/G/∞* queue system the customers arrive according to a Poisson process at rate $\lambda$, receive a service which time is a positive random variable with distribution function $G(.)$ and mean $\alpha$ and, when they arrive, find immediately an available server. Each customer service is independent from the other customers' services and from the arrivals process. The traffic intensity is $\rho = \lambda\alpha$.

A suggestion to obtain approximate results for these systems, when the exact ones are still not known, is to use a Markov renewal process, see (1-2).

Along this work, some of the approximations so obtained are evaluated.

## 2. Sojourn Time Mean Value in State *k*

For the process referred above, the sojourn time mean value in state[1] $k, k = 0,1, ...$ is given by

$$m_k = \int_0^\infty e^{-\lambda t} \left[\frac{\int_t^\infty [1 - G(x)]dx}{\alpha}\right]^k dt, k = 0,1, ... \qquad (2.1).$$

**Proposition 2.1**

$m_0 = \frac{1}{\lambda}$   (2.2).

**Obs.:** The sojourn time mean value in state 0, does not depend neither on $G(.)$ nor on the arrivals process. It depends only on the arrivals process.

---

[1] The state of the *M/G/∞* queue in instant *t* is the number of customers being served in the system at instant *t*.

**Proposition 2.2**

$$m_k \leq \frac{1}{\lambda}, \quad k = 0,1,\ldots \qquad (2.3).$$

**Dem:** It is enough to note that $\alpha^{-1} \int_t^\infty [1 - G(x)]dx \leq 1.$ □

**Obs:** So, the sojourn time mean value in any state does not exceed the one of the state 0. Put

$$E_0 = \frac{1}{\lambda}(2.4).$$

**Proposition 2.3**

$$m_k \leq \alpha \sqrt{\frac{\gamma_s^2 + 1}{2\rho(2k+1)}}, \quad k = 1,2,\ldots \qquad (2.5)$$

being $\gamma_s$ the service coefficient of variation.

**Dem:** Using the Schwartz's inequality, $m_k^2 \leq \int_0^\infty e^{-2\lambda t} dt \int_0^\infty \left[\frac{\int_t^\infty [1-G(x)]dx}{\alpha}\right]^{2k} dt =$

$\frac{1}{2\lambda\alpha^{2k}} \int_0^\infty \left[\int_t^\infty [1-G(x)]dx\right]^{2k} dt = \frac{1}{2\lambda\alpha^{2k}} \frac{2k\alpha^2}{2}(\gamma_s^2+1)\frac{\alpha^{2k-1}b_{2k-1}}{2k(2k+1)} \leq \alpha \frac{\gamma_s^2+1}{2\lambda(2k+1)}$

since, see (3),

$$\int_0^\infty \left[\int_t^\infty [1-G(x)]dx\right]^n dt = \frac{n\alpha^2}{2}(\gamma_s^2+1)\frac{\alpha^{n-1}b_{n-1}}{n(n+1)} \text{ with } b_n \leq 2,$$
$$n = 0,1,\ldots \qquad (2.6). \square$$

**Obs:** Define

$$E_1 = \alpha\sqrt{\frac{\gamma_s^2 + 1}{2\rho(2k+1)}} \quad (2.7).$$

**Proposition 2.4**

If $k \geq \frac{1}{4}\rho(\gamma_s^2 + 1) - \frac{1}{2}, E_1 \leq E_0.$

**Dem:**

$$\alpha\sqrt{\frac{\gamma_s^2+1}{2\rho(2k+1)}} \leq \frac{1}{\lambda} \Leftrightarrow \frac{\gamma_s^2+1}{2\rho(2k+1)} \leq \frac{1}{\rho^2} \Leftrightarrow \rho(\gamma_s^2+1) \leq 4k+2 \Leftrightarrow k$$
$$\geq \frac{1}{4}\rho(\gamma_s^2+1) - \frac{1}{2}. \square$$

**Proposition 2.5**

$$m_k \leq \alpha\frac{\gamma_s^2+1}{k+1}, k = 1,2,\ldots \quad (2.8).$$

**Dem:**

$$m_k \leq \int_0^\infty \left[\frac{\int_t^\infty [1-G(x)]dx}{\alpha}\right]^k dt \leq \frac{1}{\alpha^k}\frac{k\alpha^2}{2}(\gamma_s^2+1)\frac{2\alpha^{k-1}}{k(k+1)} = \alpha\frac{\gamma_s^2+1}{k+1}$$

after (2.6). □

**Obs:** Define

$$E_2 = \alpha\frac{\gamma_s^2+1}{k+1} (2.9).$$

**Proposition 2.6**

If $k \geq \rho(\gamma_s^2+1) - 1, E_2 \leq E_0$.

**Dem:**

$$\alpha\frac{\gamma_s^2+1}{k+1} \leq \frac{1}{\lambda} \Leftrightarrow \rho(\gamma_s^2+1) \leq k+1 \Leftrightarrow k \geq \rho(\gamma_s^2+1) - 1. \square$$

**Proposition 2.7**

If $k \leq 2\rho(\gamma_s^2+1) - 1, E_1 \leq E_2$.

**Dem:**

$$\frac{E_1}{E_2} = \frac{\alpha\sqrt{\frac{\gamma_s^2+1}{2\rho(2k+1)}}}{\alpha\frac{\gamma_s^2+1}{k+1}} = \frac{\sqrt{\gamma_s^2+1}(k+1)}{(\gamma_s^2+1)\sqrt{2\rho(2k+1)}} = \frac{k+1}{\sqrt{\gamma_s^2+1}\sqrt{2\rho}\sqrt{2k+1}}$$

$$\leq \sqrt{\frac{k+1}{2\rho(\gamma_s^2+1)}}; \frac{k+1}{2\rho(\gamma_s^2+1)} \leq 1 \Leftrightarrow k+1 \leq 2\rho(\gamma_s^2+1) \Leftrightarrow k \leq 2\rho(\gamma_s^2+1) - 1. \square$$

**Proposition 2.8**

If $k \geq 4\rho(\gamma_s^2+1) - 1, E_2 \leq E_1$.

**Dem:**

$$\frac{E_1}{E_2} = \frac{k+1}{\sqrt{2\rho(\gamma_s^2+1)}\sqrt{2k+1}} \geq \frac{k+1}{\sqrt{4\rho(\gamma_s^2+1)}\sqrt{k+1}} = \sqrt{\frac{k+1}{4\rho(\gamma_s^2+1)}}; \frac{k+1}{4\rho(\gamma_s^2+1)}$$

$$\geq 1 \Leftrightarrow k \geq 4\rho(\gamma_s^2+1) - 1. \square$$

The Propositions 2.4, 2.6, 2.7 and 2.8 lead to the following upper bounds choice for $m_k, k = 1,2,\ldots$:

**A)** $\rho(\gamma_s^2 + 1) > \frac{2}{3}$

$$k < \frac{1}{4}\rho(\gamma_s^2 + 1) - \frac{1}{2} \qquad m_k \leq \frac{1}{\lambda}$$

$$\frac{1}{4}\rho(\gamma_s^2 + 1) - \frac{1}{2} \leq k \leq 2\rho(\gamma_s^2 + 1) - 1 \qquad m_k \leq \alpha\sqrt{\frac{\gamma_s^2 + 1}{2\rho(2k + 1)}}$$

$$2\rho(\gamma_s^2 + 1) - 1 < k < 4\rho(\gamma_s^2 + 1) - 1 \qquad m_k \leq \min\left\{\alpha\sqrt{\frac{\gamma_s^2 + 1}{2\rho(2k + 1)}}, \alpha\frac{\gamma_s^2 + 1}{k + 1}\right\}$$

$$k \geq 4\rho(\gamma_s^2 + 1) - 1 \qquad m_k \leq \alpha\frac{\gamma_s^2 + 1}{k + 1}$$

**B)** $\frac{1}{2} < \rho(\gamma_s^2 + 1) \leq \frac{2}{3}$

$$k = 1 \qquad m_1 \leq \min\left\{\alpha\sqrt{\frac{\gamma_s^2 + 1}{6\rho}}, \alpha\frac{\gamma_s^2 + 1}{2}\right\}$$

$$k = 2,3,\ldots \qquad m_k \leq \alpha\frac{\gamma_s^2 + 1}{k + 1}$$

**C)** $\rho(\gamma_s^2 + 1) \leq \frac{1}{2}$

$$m_k \leq \alpha\frac{\gamma_s^2 + 1}{k + 1}, k = 1,2,\ldots$$

**Proposition 2.9**

If the service time distribution is *NBUE* (New Better than Used in Expectation)

$$m_k \leq \frac{\alpha}{k + \rho}, k = 1,2,\ldots \qquad (2.10)$$

**Dem:** It is enough to note that if the service time is *NBUE* with mean $\alpha$, $\int_b^\infty [1 - G(x)]dx \leq \int_b^\infty e^{-\frac{x}{\alpha}}dx$, for any $b \geq 0$. □

**Obs:** If the service time is *NWUE* (New Worse than Used in Expectation) with mean $\alpha$, $\int_b^\infty [1 - G(x)]dx \geq \int_b^\infty e^{-\frac{x}{\alpha}}dx$, for any $b \geq 0$ and

$$m_k \geq \frac{\alpha}{k + \rho}, k = 1,2,\ldots \qquad (2.11).$$

**Proposition 2.10**

If the service time distribution is *IMRL*

$$m_k \geq e^{k\left(1-\frac{2\alpha}{3\mu_2^2}\mu_3\right)} \frac{\mu_2}{\mu_2\lambda + 2k\alpha}, k = 1,2,\ldots \quad (2.12)$$

being $\mu_2$ and $\mu_3$ the 2$^{nd}$ and the 3$^{rd}$ $G(.)$ moments around the origin

**Dem:** If the service time[2] is *IMRL* (Increasing Mean Residual Life)

$$1 - G^*(x) = 1 - \frac{1}{\alpha}\int_0^x[1-G(y)]dy = \frac{\int_x^\infty[1-G(y)]dy}{\alpha} \geq e^{-\frac{2\alpha}{\mu_2}x - \frac{2\alpha}{3\mu_2^2}\mu_3 + 1}.$$

So, $\quad m_k \geq \int_0^\infty e^{-\lambda t}\left(e^{-\frac{2\alpha}{\mu_2}t - \frac{2\alpha}{3\mu_2^2}\mu_3 + 1}\right)^k dt =$

$$e^{k\left(1-\frac{2\alpha}{3\mu_2^2}\mu_3\right)}\int_0^\infty e^{-\left(\lambda + k\frac{2\alpha}{\mu_2}\right)t}dt = e^{k\left(1-\frac{2\alpha}{3\mu_2^2}\mu_3\right)}\frac{-1}{\lambda + k\frac{2\alpha}{\mu_2}}\left[e^{-\left(\lambda + k\frac{2\alpha}{\mu_2}\right)t}\right]_0^\infty$$

$$= e^{k\left(1-\frac{2\alpha}{3\mu_2^2}\mu_3\right)} \cdot \frac{\mu_2}{\mu_2\lambda + 2k\alpha}.\square$$

**Proposition 2.11**

If the service time distribution is *DFR*[3] (Decreasing Failure Rate)

$$m_k \geq e^{k\left(\frac{1-\gamma_s^2}{2}\right)}\frac{\alpha}{k+\rho}, k = 1,2,\ldots \quad (2.13).$$

**Dem:**

If the service is *DFR* $1 - G(x) \geq e^{-\frac{x}{\alpha} - \frac{\gamma_s^2}{2} + \frac{1}{2}}$.

So,

$$m_k \geq \frac{1}{\alpha^k}\int_0^\infty e^{-\lambda t}\left[\int_t^\infty e^{-\frac{x}{\alpha} - \frac{\gamma_s^2}{2} + \frac{1}{2}}dx\right]^k dt = \frac{e^{k\left(\frac{1-\gamma_s^2}{2}\right)}}{\alpha^k}\int_0^\infty e^{-\lambda t}\left[\int_t^\infty e^{-\frac{x}{\alpha}}dx\right]^k =$$

$$e^{k\left(\frac{1-\gamma_s^2}{2}\right)}\frac{\alpha}{k+\rho}.\square$$

---

[2] $G^*(x) = \frac{1}{\alpha}\int_0^x[1-G(y)]dy$ is the service time equilibrium distribution.

[3] For more details about *NBUE* (New Better than Used in Expectation), *NWUE* (New Worse than Used in Expectation), *IMRL* (Increasing Mean Residual Life) and *DFR* (Decreasing Failure Rate) distributions, important in reliability theory, see (4).

**Note:** It is not known an expression to the sojourn time value in state $k$ for the $M/G/\infty$ queue systems, except for

a) $k = 0$, for every $G(.)$, being

$$m_0 = \frac{1}{\lambda} \quad (2.14)$$

b) Every $k$, *for* exponential service time, where

$$m_k = \frac{\alpha}{k+\rho}, k = 0,1,\ldots \quad (2.15).$$

In the same circumstances, the Markov renewal process supplies the same results: in fact (2.14) is equal to (2.2) and if $G(x) = 1 - e^{-\frac{x}{\alpha}}, x \geq 0$ in (2.1) it is obtained (2.15).
-The bounds given by (2.10), (2.11), match the exact value given by (2.15). The expression (2.13) is coincident with (2.15) for $\gamma_s = 1$.

## 3. State 0 Recurrence Mean Time

For the Markov renewal process, the state 0 mean recurrence time[4] is given by

$$\mu_0 = \frac{1}{\lambda}\left[1 + \sum_{j=1}^{\infty} \prod_{k=1}^{j} \frac{\lambda m_k}{1 - \lambda m_k}\right] \quad (3.1).$$

**Proposition 3.1**

If $\rho \leq \frac{1}{\gamma_s^2+1}$, $\mu_0 \leq \frac{e^{\rho(\gamma_s^2+1)}}{\lambda} \quad (3.2)$

**Dem:** To use an upper bound of $m_k$ in (3.1) it is necessary to certify that it is lesser than $\frac{1}{\lambda}$. The condition $\rho(\gamma_s^2 + 1) \leq 1$, due to Proposition 2.6, guaranties that $E_2$ fulfills that request for $k \geq 1$.

So $\mu_0 \leq \frac{1}{\lambda}\left[1 + \sum_{j=1}^{\infty} \prod_{k=1}^{j} \frac{\frac{\rho(\gamma_s^2+1)}{k+1}}{1 - \frac{\rho(\gamma_s^2+1)}{k+1}}\right] = \frac{1}{\lambda}\left[1 + \sum_{j=1}^{\infty} \prod_{k=1}^{j} \frac{\rho(\gamma_s^2+1)}{k+1-\rho(\gamma_s^2+1)}\right] \leq \frac{1}{\lambda}\left[1 + \sum_{j=1}^{\infty} \frac{[\rho(\gamma_s^2+1)]^j}{j!}\right] = \frac{e^{\rho(\gamma_s^2+1)}}{\lambda}$. □

**Obs:** For the $M/G/\infty$ queue systems

$$\mu_0 = \frac{e^\rho}{\lambda} \quad (3.3).$$

So, in these conditions, the relative error arising from considering (3.2) instead of (3.1) is

$$\frac{e^{\rho(\gamma_s^2+1)}}{\lambda} - \frac{e^\rho}{\lambda}}{\frac{e^\rho}{\lambda}} = e^{\rho\gamma_s^2} - 1 \leq e^{\frac{\gamma_s^2}{\gamma_s^2+1}} - 1 < e - 1.$$

---

[4] It is in fact the $M/G/\infty$ queue busy cycle mean time, see (5).

But $e^{\frac{\gamma_s^2}{\gamma_s^2+1}} - 1 \leq r \Leftrightarrow \frac{\gamma_s^2}{\gamma_s^2+1} \leq \log(r+1) \Leftrightarrow \gamma_s^2 \leq \frac{\log(r+1)}{1-\log(r+1)}$.

That is: if $\rho(\gamma_s^2 + 1) \leq 1$, the relative error arising from taking the bound given by (3.2) instead of the true value given by (3.3) for $\mu_0$ is such that:

a) $\varepsilon \leq e^{\frac{\gamma_s^2}{\gamma_s^2+1}} - 1$,

b) $\varepsilon = 0$ if $\gamma_s^2 = 0$,

c) $\varepsilon < e - 1$,

d) $\varepsilon \leq r$ $(r < e - 1)$ since $\gamma_s^2 \leq \frac{\log(r+1)}{1-\log(r+1)}$.

So, requesting that $\varepsilon$ is lesser than a given r, it results a criterion to measure the goodness of the $m_k$ approximation by $E_2$ for a certain $\gamma_s^2$, since $\rho(\gamma_s^2 + 1) \leq 1$.

Being B the $M/G/\infty$ queue busy period length, see (5-6),

$$E[B] = \frac{e^\rho - 1}{\lambda} \quad (3.4).$$

For the Markov renewal process, since $\rho(\gamma_s^2 + 1) \leq 1$,

$$E[B] = \frac{e^{\rho(\gamma_s^2+1)}-1}{\lambda} \quad (3.5).$$

Now, the relative error own to take (3.5) instead (3.4), is

$$\frac{\frac{e^{\rho(\gamma_s^2+1)}-1}{\lambda} - \frac{e^\rho-1}{\lambda}}{\frac{e^\rho-1}{\lambda}} = \frac{e^{\rho(\gamma_s^2+1)}-e^\rho}{e^\rho-1} = \frac{e^{\rho\gamma_s^2}-1}{1-e^{-\rho}} \leq \frac{e^{\frac{\gamma_s^2}{\gamma_s^2+1}}-1}{1-e^{-\rho}} < \frac{e-1}{1-e^{-\rho}}. \text{ So}$$

a) $\delta \leq \frac{e^{\frac{\gamma_s^2}{\gamma_s^2+1}}-1}{1-e^{-\rho}}$,

b) $\varepsilon = 0$ if $\gamma_s^2 = 0$,

c) $\delta < \frac{e-1}{1-e^{-\rho}}$,

d) $\delta \leq r \left( r < \frac{e-1}{1-e^{-\rho}} \right)$ since $\gamma_s^2 \leq \frac{\log(r(1-e^{-\rho})+1)}{1-\log(r(1-e^{-\rho})+1)}$.

So, the bounds for $\delta$ are greater than the obtained to $\varepsilon$. Then it is preferable to use a criterion based on $\varepsilon$ than on $\delta$, to measure the goodness of the approximation of $m_k$ by $E_2$.

## 4. Mean Number of Entries in State $k$ between Two Entries in State 0

For the Markov renewal process, the mean number of entries in state $k$ between two entries in state 0 is

$$v_k = \lambda^{k-1} \frac{m_1 \ldots m_k}{(1-m_1)\ldots(1-m_k)}, k = 1,2,\ldots \qquad (4.1).$$

**Proposition 4.1**

If $\rho \leq \frac{1}{\gamma_s^2+1}$

$$v_k \leq (k+1) \frac{\rho^{k-1}(\gamma_s^2+1)^{k-1}}{k!}, k = 1,2,\ldots \qquad (4.2).$$

**Obs.:** Values for $v_k, k = 1,2,\ldots$ for the $M/G/\infty$ queue system is not known,

From (2.2), (2.9) and (4.2) it follows that

$$m_k v_k \leq \frac{\alpha \rho^{k-1}(\gamma_s^2+1)^k}{k!}, k = 0,1,\ldots \quad (4.3)$$

Since $\rho(\gamma_s^2+1) \leq 1$.

For the $M/G/\infty$ queue system

$$m_k v_k = \frac{\alpha \rho^{k-1}}{k!}, k = 0,1,\ldots \qquad (4.4).$$

But $\dfrac{\frac{\alpha\rho^{k-1}(\gamma_s^2+1)^k}{k!} - \frac{\alpha\rho^{k-1}}{k!}}{\frac{\alpha\rho^{k-1}}{k!}} = (\gamma_s^2+1)^k - 1$, that is null for $\gamma_s = 0$ or $k = 0$ and

increases with $k$ if $\gamma_s^2 > 0$.

**Note:** For $\rho \leq 1$ *and* $\gamma_s^2 = 0$ the Markov renewal process supplies the following results:

a) $m_0 = \frac{1}{\lambda}$,

b) $m_k \leq \frac{\alpha}{k+1}, k = 1,2,\ldots$

c) $\mu_0 \leq \frac{e^\rho}{\lambda}$,

d) $E[B] \leq \frac{e^{\rho}-1}{\lambda}$,

e) $v_k \leq (k+1)\frac{\rho^{k-1}}{k!}$, $k = 1,2,...$

f) $m_k v_k = \frac{\alpha \rho^{k-1}}{k!}$, $k = 0,1,...$

So, the value obtained for $m_0$ coincides with the $M/G/\infty$ one. And the bounds obtained for $\mu_0$, $E[B]$ and $m_k v_k$ coincide with the true values for the same $M/G/\infty$ quantities. But the bounds obtained for $m_k$ and $v_k$ coincide with the true value obtained when the service time distribution is exponential, and the traffic intensity is 1. In opposition, the bound got for $m_k$ cannot coincide with the one given by (2.15) for $\rho < 1$. So, it is excluded the hypothesis of having an expression for $m_k$ independent from the service time distribution and equal to the one given by (2.15). Then, only rarely the Markov renewal process gives values for $\mu_0$, $E[B]$ and $m_k v_k$ identical to the $M/G/\infty$ ones.

If $\rho(\gamma_s^2 + 1) \leq 1$ it is possible, after the Markov renewal process, to get upper bounds for the $M/G/\infty$ system quantities $\mu_0$, $E[B]$ and $m_k v_k$. So, it is admissible to consider that at least $E_2$, beyond being a $m_k$ upper bound for the Markov renewal process, also plays the same role for the $M/G/\infty$ queue.

Note still that if $\gamma_s^2 = 0$, for the Markov renewal process

$$m_k = \int_0^\alpha e^{-\lambda t}\left[1-\frac{t}{\alpha}\right]^k dt, \quad k = 0,1 ... \qquad (4.5).$$

So, for $k \geq 1$, $m_k \leq \int_0^\alpha \left(1-\frac{t}{\alpha}\right)^k dt = \left[\frac{-\alpha}{k+1}\left(1-\frac{t}{\alpha}\right)^{k+1}\right]_0^\alpha$,

That is $m_k \leq \frac{\alpha}{k+1}$, $k = 1,2,...$.

But, requesting that $\frac{\alpha}{k+1} \leq \frac{1}{\lambda} \Leftrightarrow k \geq \rho - 1$ that leads to $\rho - 1 \leq 0 \Leftrightarrow \rho \leq 1$.

$$m_1 = \int_0^\alpha e^{-\lambda t}\left(1-\frac{t}{\alpha}\right) dt = \left[-\frac{e^{-\lambda t}}{\lambda}\left(1-\frac{t}{\alpha}\right)\right]_0^\alpha -$$

$\int_0^\alpha -\frac{e^{-\lambda t}}{\lambda}\left(-\frac{1}{\alpha}\right) dt = \frac{1}{\lambda} - \frac{1}{\rho}\int_0^\alpha e^{-\lambda t} dt = \frac{1}{\lambda} - \frac{1}{\rho}\left[\frac{e^{-\lambda t}}{-\lambda}\right]_0^\alpha = \frac{1}{\lambda} - \frac{1}{\rho}\left(-\frac{e^{-\rho}}{\lambda} + \frac{1}{\lambda}\right)$. So,

$$m_1 = \alpha \frac{\rho + e^{-\rho} - 1}{\rho^2} \qquad (4.6).$$

And, integrating by parts,

$$m_{k+1} = \frac{1}{\lambda} - \frac{k+1}{\rho} m_k, \text{ k=1, 2, ... (4.7).}$$

With (4.6) and (4.7) it is possible to obtain $m_k$, $k = 1,2,\ldots$ for $\gamma_s^2 = 0$ and it is possible to conclude that, in this case, (2.15) does not hold.

**Proposition 4.2**

If the service time distribution is *NBUE*

a) $\mu_0 \leq \frac{e^\rho}{\lambda}$,

b) $E[B] \leq \frac{e^\rho - 1}{\lambda}$,

c) $v_k \leq (k+1)\frac{\rho^{k-1}}{k!}$, $k = 1,2,\ldots$

d) $m_k v_k \leq \frac{\alpha \rho^{k-1}}{k!}$, $k = 0,1,\ldots$

**Obs:** The bounds obtained for $\mu_0$, $E[B]$ and $m_k v_k$ coincide with the true value of these quantities for the $M/G/\infty$ queue.

If the service time distribution is *NWUE*

a) $\mu_0 \geq \frac{e^\rho}{\lambda}$,

b) $E[B] \geq \frac{e^\rho - 1}{\lambda}$,

c) $v_k \geq (k+1)\frac{\rho^{k-1}}{k!}$, $k = 1,2,\ldots$

d) $m_k v_k \geq \frac{\alpha \rho^{k-1}}{k!}$, $k = 0,1,\ldots$

with a comment identical to the one in the case *NBUE*.

So, it is admissible that $\frac{\alpha}{k+\rho}$, $k = 0,1,\ldots$ is an upper bound (lower bound) for the true value of $m_k$ in the $M/G/\infty$ queue systems in the case of *NBUE* (*NWUE*) service time distributions.

After (2.1) and integrating by parts

$m_{k+1} = \int_0^\infty e^{-\lambda t} \left[\frac{\int_t^\infty [1-G(x)]dx}{\alpha}\right]^{k+1} dt = \left[-\frac{e^{-\lambda t}}{\lambda}\left(\frac{\int_t^\infty [1-G(x)]dx}{\alpha}\right)^{k+1}\right]_0^\infty -$

$\int_0^\infty -\frac{e^{-\lambda t}}{\lambda}(k+1)\left[\frac{\int_t^\infty [1-G(x)]dx}{\alpha}\right]^k \frac{G(t)-1}{\alpha} dt = \frac{1}{\lambda} - \frac{k+1}{\rho}\int_0^\infty e^{-\lambda t}\left[\frac{\int_t^\infty [1-G(x)]dx}{\alpha}\right]^k (1-G(t))dt \geq \frac{1}{\lambda} - \frac{k+1}{\rho} m_k$. So

$$m_{k+1} \geq \frac{1}{\lambda} - \frac{k+1}{\rho} m_k, \quad k = 1,2,\ldots \quad (4.8).$$

**Note:** According to (4.7), when the service is constant, the equality holds in (4.8).

## 5. Sojourn Time in State k Distribution

The sojourn time in sate k distribution function for the Markov renewal process is

$$C_k(t) = 1 - e^{-\lambda t}\left[\frac{\int_t^\infty [1-G(x)]dx}{\alpha}\right]^k, t \geq 0, k = 0,1,\ldots \quad (5.1).$$

Evidently,

**Proposition 5.1**

$$C_k(t) \geq 1 - e^{-\lambda t}, t \geq 0, k = 0,1,\ldots \quad (5.2).$$

**Proposition 5.2**

If the service time distribution is exponential

$$C_k(t) = 1 - e^{-\frac{1}{\alpha}(k+\rho)t}, t \geq 0, k = 0,1,\ldots \quad (5.3).$$

**Obs.:** This result is coincident with the one known for the $M/G/\infty$ queue.

**Proposition 5.3**

$$C_0(t) = 1 - e^{-\lambda t}, t \geq 0, \quad (5.4).$$

**Obs.:** Result obvious for any $M/G/\infty$ queue and for any queue with Poisson arrivals process.

**Proposition 5.4**

If the service time distribution is *NBUE*

$$C_k(t) \geq 1 - e^{-\frac{1}{\alpha}(k+\rho)t}, t \geq 0, k = 0,1,\ldots \quad (5.5).$$

**Obs.:** As emphasized before, (5.5), beyond supplying a lower bound for $C_k(t)$ in the Markov renewal process, also gives a lower bound for that quantity in the $M/G/\infty$ system for the case of *NBUE* service time.

If the service time distribution is *NWUE*

$$C_k(t) \leq 1 - e^{-\frac{1}{\alpha}(k+\rho)t}, t \geq 0, k = 0,1,\ldots \quad (5.6)$$

And it is pertinent a comment analogous to the former one with the change of lower bound by upper bound.

**Proposition 5.5**

If the service time distribution is *IMRL*

$$C_k(t) \leq 1 - e^{-\lambda t + k\left(-\frac{2\alpha}{\mu_2}t - \frac{2}{3}\frac{\alpha}{\mu_2^2}\mu_3 + 1\right)}, t \geq 0, k = 0,1,\ldots \quad (5.6)$$

**Proposition 5.6**

If the service time distribution is *DFR*

$$C_k(t) \leq 1 - e^{-\frac{1}{\alpha}(k+\rho)t + k\frac{1-\gamma_s^2}{2}}, t \geq 0, k = 0,1,\ldots \quad (5.7).$$

## Conclusions

When analytical exact results are not available, numerical methods are used to try to find approximations for the interesting quantities under study. It is what is done in this work for the *M/G/∞* queue, trying to approximate it for a Markov renewal process. An alternative is using simulation methods. For this approach see, for instance (6, 7).

Still another is to determine service time distributions for which it is possible to determine most of the interesting quantities for the *M/G/∞* queue. This is made solving differential equations induced for the study of the transient behavior, see (8-10).

## References

[1]  M.F. Ramalhoto and D.H. Girmes, Markov Renewal Approach to Counter Theory, in: J.P. Barra et al (Eds.), Recent Developments in Statistics, North Holland, 581-590, 1977.
[2]  M.F. Ramalhoto, Breves Considerações Sobre Aproximações e Limites em Sistemas do Tipo *GI/G/∞*.IST, CEAUL, Lisboa, 1985.
[3]  Y.S. Sathe, Improved bounds for the variance of the busy period of the *M/G/∞* queue, Advances in Applied Probability, 17, 913-914, http://dx.doi.org/10.2307/1427096. 1985.
[4]  S. Ross, Stochastic Processes, Wiley, New York, 1983.
[5]  M.A.M. Ferreira and M. Andrade, Busy period and busy cycle distributions and parameters for a particular *M/G/∞* queue system, American Journal of Mathematics and Statistics, 2(2) (2012), 10-15, http://article.sapub.org/10.5923.j.ajms.2012 0202.03.html.
[6]  M.A.M. Ferreira, Simulação computacional de sistemas com infinitos servidores, Revista de Estatística, 3(3) (1998), 5-28.
[7]  M.A.M. Ferreira and M. Andrade, Infinite servers queue systems computational simulation, 12[th] Conference on Applied Mathematics- APLIMAT 2013, Bratislava, (2013).
[8]  M.A.M. Ferreira and M. Andrade, The ties between the *M/G/∞* queue system transient behavior and the busy period, International Journal of Academic Research, 1(1) (2009), 84-92.
[9]  M.A.M. Ferreira, Differential equations important in the *M/G/∞* queue system transient behaviour busy period study, International Conference in Applied Mathematics- APLIMAT 2005, Proceedings, Bratislava, (2005), 119-132.

[10]  M.A.M. Ferreira and M. Andrade, $M/G/\infty$ system transient behavior with time origin at the beginning of a busy period mean and variance, APLIMAT- Journal of Applied Mathematics, 3(3) (2010), 213-221.

[11]  M.A.M. Ferreira, Aplicação da equação de Riccati ao estudo do período de ocupação do sistema $M/G/\infty$, Revista de Estatística, 3(1) (1998), 23-28.

[12]  M.A.M. Ferreira and M. Andrade, $M/G/\infty$ queue system parameters for a particular collection of service time distributions, AJMCSR-African Journal of Mathematics and Computer Science Research, 2(7) (2009), 138-141.

[13]  M.A.M. Ferreira, M. Andrade, and J.A. Filipe, The Riccati equation in the $M/G/\infty$ system busy cycle study, Journal of Mathematics, Statistics and Allied Fields, 2(1) (2008).

[14]  M.A.M. Ferreira, M. Andrade and J.A. Filipe, Networks of queues with infinite servers in each node applied to the management of a two echelons repair system, China-USA Business Review, 8(8) (2009), 39-45.

[15]  M.J. Carrillo, Extensions of Palm's theorem: A review, Management Science, 37(6) (1991), 739-744.